\documentclass[12pt]{amsart}
\usepackage{geometry}                
\usepackage{graphicx}
\usepackage{amssymb}
\usepackage{color}
\usepackage{pstricks}

\newtheorem{theorem}{Theorem}
\newtheorem{corollary}[theorem]{Corollary}    
\newtheorem{definition}{Definition}
\newtheorem{proposition}{Proposition}
\newtheorem{lemma}{Lemma}
\newtheorem{remark}{Remark}

\usepackage{graphicx, color}
\usepackage{amssymb}
\usepackage{epstopdf}
\title{Permanents and Determinants, Weighted Isobaric Polynomials and Integral Sequences}
\author{Huilan Li and Trueman MacHenry}

\begin{document}

\begin{abstract}
In this paper we construct two types of Hessenberg matrices with the properties that every weighted isobaric polynomial (WIP) appears as a determinant of one of them,  and as the permanent of the other.  Every integer sequence which is linearly recurrent is representable by (an evaluation  of) some linearly recurrent sequence of WIPs.  WIPs are symmetric polynomials written on the elementary symmetric polynomial basis.  Among them are the generalized Fibonacci polynomials and the generalized Lucas polynomials,  which already have these sweeping representing properties.  Among the integer sequences discussed are the Chebychev polynomials of the 2nd kind, the Stirling numbers of the 1st and 2nd kind, the Catalan numbers, and the triangular numbers,  as well as all sequences which are either multiplicative arithmetic functions or additive arithmetic functions.
\end{abstract}
\vspace{0.5cm}

\maketitle



\section{Introduction}
In this paper,  we give a general method for finding permanental and determinantal representations of many families of integer sequences These include all of the integer sequences in recent papers of Kirgisiz and Sahin, for example, \cite{KS}. However,  the methods used here are general and are part of an overarching theoretical structure.  The tools that we use involve the ring of symmetric polynomials and the methods are ones that should be congenial to workers in algebraic combinatorics,

The term \textit{permanent} dates from Cauchy, 1821,  and the modern usage of the term as permanent of a matrix goes back to Muir in 1882 \cite{Mu}. A permanent is computed very much like a determinant except that one ignores the parity of the elements acting from  the symmetric group.  The interest in the permanent of a matrix is at least  two-fold.  On the one hand, it has applications in graph theory,  on the other hand,  it is used in quantum physics.  Since its applications are computational, finding economical computing techniques is a desirable end.  Clearly the more zeros that occur in a matrix,  the easier will be the computation of its permanent as well as its determinant.

\textit{Hessenberg} Matrices, \textit{upper} and \textit{lower}, are matrices that are nearly triangular.  Specifically,  an \textit{upper Hessenberg} matrix is one with a not-necessarily zero sub-main diagonal,  otherwise there are only zeros below the main diagonal; a \textit{lower Hessenberg} matrix has a not-ncessarily zero super-main diagonal,  otherwise there are only zeros above the main diagonal.

The representation theorems in \cite{KS} are of this sort,  as will be ours as well. After having described a very general class of polynomials (\textit{weighted isobaric} polynomials, WIPS) whose evaluations produce a large number of interesting integer sequences (section 3) simply by varying the parameters in these polynomials,  we then construct (section 8) two special Hessenberg matrices, one of which has these polynomials as permanents,  the other has these polynomials as determinants.

The paper begins  with a presentation of the theory of isobaric polynomials (section 2)  following the earlier papers of Li, MacHenry, Tudose and Wong in \cite{TM,TM1,MT,MT2,MW,MW2,LM}.  Isobaric polynomials are symmetric polynomials written on the elementary symmetric polynomial  basis. Two especially important sequences of isobaric polynomials are the \textit{Generalized Fibonacci Polynomial} sequence (GFP) and the \textit{Generalized Lucas Polynomial} sequence (GLP) introduced in \cite{TM}.  (Both of these terms, ``generalized Fibonacci'' and ``generalized Lucas''  have been used for many years, and continue to be used,  but they refer to two-variable polynomials,  while our polynomials have arbitrarily many variables and include the two-variable case) and appeared for the first time in \cite{TM, TM1}.   GFPs and GLPS are examples of WIPs, and both sequences are bases for the ring of symmetric polynomials.

The theory of isobaric polynomials can be thought of as a collection of \textit{packages} of results following from the specification of a monic polynomial: --- for the purposes of this paper ---  with integer coefficients.  We call  this polynomial the \textit{core polynomial}. We distinguish between a \textit{generic }monic polynomial, one with variable coefficients ,  and a \textit{numerical} polynomial in which the coefficients have been evaluated over, say, the ring of integers. A package consists of a \textit{generic  core} polynomial (GCP) and its companion matrix (GCM), and of an extension of the  GCM,  the i\textit{nfinite companion matrix} (ICM), whose right hand column contains the positively and negatively indexed GFPs,  and its diagonal elements give the GLPs \cite{MT}.  Each evaluation of the GCM induces the generic linear recursion of degree $n$ (GLR(n)) \cite{MW};  it also induces a multiplicative arithmetic function (MF),  and an additive arithmetic function  (AddF) in the ring (unique factorization domain) of arithmetic functions \cite{CE,MW2,LM},  and the GCP gives families of MF and AddF,  \cite{MW2}. Moreover ,  we can also think of the rational extension ring (field,  if CP is irreducible) obtained from the GCP or the CP as being part of this package, so, in particular, when CP is irreducible,  a class of number fields (NF(n)) is induced \cite{MW}.  The Frobenius character theorem can be represented in terms of WIPs, so the GCP of degree $n$ induces the character table of $Sym(n)$ (Ch(n)), \cite{LM} .  Thus given a GCP of degree $n$, we can think of the package as the collection, $$ (GCP, GCM, ICM, GLR(n), GFP(n), GLP(n), MF(n), AddF(n), NF(n), ChT(n) ).$$  This can be thought of as the theory arising from the classical result that the coeficients of a monic polynomial are elementary symmetric functions (ESP) of its roots, that is an answer to the question, what follows from the way the roots determine the coefficients of a monic polynomial? This is somewhat analogous to the way in which Galois theory arises from asking the question:  in what way do the coefficients determine the roots of the polynomial?

For the history of the term \textit{isobaric} itself,  see  Read, \cite{Read}, Redfield, \cite{Red} and P\'{o}lya,\cite{P}.

  \vspace{0.5cm}

 \noindent Section 1.  Introduction\\
Section 2. \textit{core polynomial}, the \textit{ companion matrix} and \textit{isobaric polynomials}.\\
 Section 3. \textit{weighted isobaric polynomials}  and the \textit{infinite companion matrix}.\\
 Section 4. \textit{Different Matrix} and related results. \\
 Section 5.  \textit{Convolution Product} .\\
 Section 6.  \textit{ Isobaric Logarithm}, \textit{Isobaric Exponential Operators} and  \textit{Isobaric\,Trigonometry}.\\
 Section 7.  \textit{Multiplicative Arithmetic Functions}.\\
 Section 8.  \textit{Hessenberg Matrices}.\\
 Section 9. \textit{Permanental} and \textit{Determinantel} Representation.\\
 Section 10. \textit{Representability of Integer Sequences}.\\
 Section 11. Examples.\\
\section{Isobaric Polynomials}

Let
$$\mathcal{C}(X) = X^k-t_1X^{k-1}-\cdots-t_1$$
be a \textit{generic} degree-k  monic polynomial, that is,  we consider the coefficients to be variables which can be evaluated over a suitable ring.  In this paper,  that ring will be the ring of integers, $\mathbb{Z}$.  We call this polynomial the \textit{core polynomial }. It is with respect to this polynomial that we shall make the following definitions.  For this polynomial we construct the \textit{companion matrix}:

$$A = \left(\begin{array}{ccccc}0 & 1 & 0 & \cdots & 0 \\0 & 0 & 1 & \cdots & 0 \\\vdots & \vdots & \vdots & \ddots & \vdots \\0 & 0 & 0 & \cdots & 1 \\t_k & t_{k-1} & t_{k-2} & \cdots & t_1\end{array}\right).$$

The Companion Matrix will be used to generate some remarkable polynomials which we now define.

 \begin{definition}\label{isoP} An isobaric polynomial is a polynomial on the variables $t_1, t_2, \ldots , t_k $ for  $k \in \{1,2,\ldots,n, \ldots \}$,  with coefficients,  for purposes of this paper, \, in $\mathbb{Z}$,  of the form

$$P_{k,n}(t_1,t_2 \ldots , t_k)=  \sum_{\alpha \vdash n} {C_\alpha }t_1^{\alpha_1} t_2^{\alpha_2} \cdots t_k^{\alpha_k} $$
where  $\alpha = \{\alpha_1, \alpha_2, \ldots ,\alpha_k\}$  and  where,  $\alpha \vdash n$ means that
$\sum_{j=1} ^k j \alpha_j = n$;  that is , in a standard  partition notation, where  $(1^{\alpha_1}, 2^{\alpha_2}, \ldots, k^{\alpha_k})$ is a partition of  $n$. Thus isobaric polynomials are polynomials indexed  by partitions of the natural numbers, or, equivalently,  by Young diagrams. $n$ is called the \textit{isobaric degree} of the polynomial.
\end{definition}

 \vspace{0.5cm}

 \begin{theorem}\label{myisoring} The isobaric polynomials form a ring which is isomorphic to the ring of symmetric polynomials.  In fact,  the isobaric polynomials are simply the symmetric polynomials written on the elementary symmetric polynomial basis.  The isomorphism is given by letting $t_j$ go to $(-1)^j j$-th elementary symmetric polynomial for each grading   $k$ of the graded ring of symmetric polynomials.

 \end{theorem}

\section{Weighted Isobaric Polynomials and the Infinite Companion Matrix}

 \begin{definition}\label{myWIP} A \textit{Weighted Isobaric Polynomial} is an isobaric polynomial given by the following explicit expression:

$$P_{\omega,k,n}(t_1^{\alpha_1} t_2^{\alpha_2} \cdots t_k^{\alpha_k})= \sum_{\alpha \vdash n} \left( \begin{matrix}
 |\alpha|  \\  \alpha_1,\ldots,\alpha_k \end{matrix} \right)  \frac{\sum_{j=1}^k \omega_j\alpha_j}{|\alpha|}t_1^{\alpha_1}\cdots t_k^{\alpha_k}$$
where   $\omega$ is the weight vector $(\omega_1,\omega_2, \ldots, \omega_n, \ldots), \omega_j \in\mathbb{Z}$  and  $|\alpha|= \alpha_1 + \cdots + \alpha_k.$  The coefficient   ${C_\alpha }$ is then just $$\sum_{\alpha \vdash n} \left( \begin{matrix}  |\alpha|  \\  \alpha_1,\ldots,\alpha_k \end{matrix} \right)  \frac{\sum_{j=1}^k \omega_j\alpha_j}{|\alpha|}.$$

 \end{definition}


Generating  functions for the weighted isobaric polynomials are given by:

$$\Omega(y) = 1 + \frac{\omega_1t_1y +\omega_2t_2y^2 + \omega_3t_3y^3 + \cdots+\omega_kt_ky^k}{1-p(y)},$$ where $p(y) = t_1y+t_2y^2+t_3y^3+\cdots+t_ky^k.$
 An easy induction gives the following very useful proposition.
   \vspace{0.5cm}

 \begin{proposition}\label{myLR} (Linear Recursion Property).  $$P_{\omega,k,n} =\sum_{j=1}^k t_j P_{\omega,k,n-j} (\mathbf{t}),$$  where    $\mathbf{t} = (t_1, \ldots, t_k)$.

 \end{proposition}

 \vspace{0.5cm}
 Two important sequences are the \textit{Generalized Fibonacci Polynomials} GFP, and the \textit{Generalized Lucas Polynomials} GLP.  The weight vector for the GFPs is $(1,1,\ldots,1,\ldots)$ and for the GLPs, $(1,2,\ldots,n,\ldots).$  Thus,

  \vspace{0.5cm}
 \begin{corollary}\label{myGFP,GJP } (Explicit formulae for GFP and GLP)

 1.   $$F_{k,n} = \sum_{\alpha \vdash n}\left( \begin{array}{c}
     |\alpha|    \\
      \alpha_1, \ldots, \alpha_k
\end{array}\right)  t_1^{\alpha_1} t_2^{\alpha_2} \cdots t_k^{\alpha_k}.$$

2.  $$G_{k,n} = \sum_{\alpha \vdash n} \frac{n}{|\alpha|} \left(\begin{array}{c}
     |\alpha|    \\
      \alpha_1, \ldots, \alpha_k
\end{array} \right)  t_1^{\alpha_1} t_2^{\alpha_2} \cdots t_k^{\alpha_k}.$$
  \vspace{0.5cm}
  \end{corollary}
 The generating function for elements in  $\mathcal{W}$ given above specializes to
   \vspace{0.5cm}

 1. $$\Omega_F(y) = \frac{1}{1-p(y)}  \quad \hbox{for }  GFP, $$
 and to

 2. $$\Omega_G(y) = \frac{1+t_2y^2+\cdots+(k-1)t_ky^k}{1-p(y)} \quad  \hbox{for }  \, GLP.$$

   \vspace{0.5cm}

Another way of producing these polynomials is through the \textit{ infinite companion }matrix, which  we now define.
  \vspace{0.5cm}

 \begin{definition}\label{my ICM } The Infinite companion matrix  $A^\infty$ is the matrix defined by allowing the companion matrix $A$ to operate on the row vectors of $A$,  that is on itself.  We clearly reproduce $A$ itself this way,  and if we sequentially adjoin the orbit vectors under this operation as new rows of  $A$,  we obtain infinitely many new rows southward.  If we also agree that  $ t_k\neq 0$,  so that  $A$ is non-singular,  and operate on the first row vector of $A$ with  $A^{-1}$,  and adjoin this orbit northward to  $A$,  we obtain infinitely many rows northward,  so that  the resulting matrix is an  $(\infty \times k)$-matrix.   As we shall see,  it is a rather remarkable matrix which contains a great deal of information. Below we give a ``picture'' of  $A^\infty$.
 \end{definition}

    $$A_k^\infty  = \left( \begin{matrix}

     \vdots &  \ddots & \vdots & \vdots \\
      (-1)^{k-1}S_{(-2,1^{(k-1)})} & \cdots  & -S_{(-2,1)} & S_{(-2)}  \\
      (-1)^{k-1}S_{(-1,1^{(k-1)})} & \cdots  & -S_{(-1,1)} & S_{(-1)}  \\
      (-1)^{k-1}S_{(0,1^{(k-1)}))}  &\cdots   &  -S_{(0,1)} & S_{(0)}  \\
      (-1)^{k-1}S_{(1,1^{(k-1)})}  & \cdots  &  -S_{(1,1)} & S_{(1)}  \\
      (-1)^{k-1}S_{(2,1^{(k-1)})}  & \cdots   &  -S_{(2,1)} & S_{(2)}  \\
      (-1)^{k-1}S_{(3,1^{(k-1)})}  & \cdots   &  -S_{(3,1)} & S_{(3)}  \\
      (-1)^{k-1}S_{(4,1^{(k-1)})} & \cdots   &  -S_{(4,1)} & S_{(4)}  \\
                    \vdots                      &  \ddots    &    \vdots     &    \vdots

   \end{matrix} \right) = ((-1)^{k-j} S_{(n,1^{k-j})}. $$

 \begin{lemma}\label{my CMLR } Let   $v$ be a row vector with $k$ components.  The orbit of $v$ under the action of  $A$ on the right, giving the ordered set  $vA^n$, is a linearly recursive sequence with recursion parameters  $\{ t_1,t_2,\ldots,t_k \} $.  In fact,  we get all k-th order linear recursions in this way, \cite{MW}.
\end{lemma}\qed

 \vspace{0.5cm}

 \begin{corollary}\label{mycolMF } Each  column of the infinite companion matrix is a linearly recursive sequence of degree  k.\end{corollary}
  \vspace{0.5cm}

 \proof It is easy to see that  the orbit of any vector under the operation of  $A$ is necessarily a linar recursion with the set $\{t_j\}$ as the recursion parameters.  \qed

 \vspace{0.5cm}

 \begin{theorem}{\cite{LM}}\label{my propsCM }
 \vspace{0.5cm}

We assume that  $t_k \neq 0$,  so that,  $A_k^{-1}$ is defined.
 \vspace{0.5cm}

 1.   The powers of $A_k$ constitute an infinite cyclic group.
  \vspace{0.5cm}

 2.  The $k \times k$ contiguous blocks of  $A_k^\infty$ are the powers  $A_k^n$ of $A_k$.  In particular,  the $k \times k$  block whose lower right-hand corner is $S_{(0)}$ is the identity matrix, and the $k \times k$ block whose lower right-hand corner is  $S_{(1)}$ is just  $A_k$.

  \vspace{0.5cm}
  3. The entries in  $A_k^\infty$ are Schur-hook polynomials with arm-length $n$ and leg-length  $r$.
   \vspace{0.5cm}

  4.  The right-hand column is the linearly recursive sequence of GFPs,  the Generalized Fibonacci Polynomials.  The traces of the infinite cyclic group generated by  $A_k$,  that is,  the diagonals of $A_k^\infty$,  is the linearly recursive sequence of GLPs, that is,  the Generalized Lucas Polynomials.
   \vspace{0.5cm}

  5.  All of the columns of  $A_k^\infty$ are linear recursions with recursion parameters $\{t_1, t_2, \ldots, t_k \}$,  and every linear recursion of degree $k$ can be obtained by a suitable choice of k and evaluation of the parameters. In particular,  when  $k=2$ and $t_1 = t_2 = 1$,  then GFP becomes the Fibonacci sequence,  and  GLP becomes the Lucas sequence.
  \vspace{0.5cm}

 6.   Letting the collection of weighted isobaric polynomials be denoted by  $\mathcal{W}$, we can define a group operation on  $\mathcal{W}$ by applying the usual component-wise addition of vectors to the  weight vectors.  Thus given weights   $\omega$ and $\omega'$,  we can define a new weight  vector $\omega + \omega'$.  It is trivial to see that this induces a group structure on $\mathcal{W}$.  Call this group  $\mathcal{W(\omega)}$. Also ring addition of two weighted isobaric polynomials with possibly different weights gives a group structure on $\mathcal{W}$.  Call it $\mathcal{W(+)}$

 \end{theorem}
  \vspace{0.5cm}

  \begin{theorem}\label{myW+ }  $\mathcal{W(\omega)} \cong \mathcal{W(+)}$.  In fact they are identical.    \qed
\end{theorem}


\section{The Different Matrix}
  \vspace{0.5cm}

Consider the derivative of the core polynomial,
$$C'(X) = kX^{k-1} - (k-1)t_1X^{k-2} - \cdots-t_{k-1}.$$
  \vspace{0.5cm}
Define the vector    $d_k = (-t_{k-1} , -2t_{k-2}, \ldots, -(k-1)t_{1}, k),$  and construct the \textit{Different} matrix

 $$D_k = \left( \begin{matrix}
d_kA^0\\ d_kA\\ \vdots\\ d_kA^{k-1}
 \end{matrix} \right),$$

\noindent the matrix obtained by operating on  $d_k$ with the companion matrix.  For example,  when  $k=3$,

$$ D_3 = \left( \begin{matrix}
-t_2 & -2t_1 & 3 \\
3t_3 & 2t_2 & t_1 \\
t_1t_3 & 3t_3+t_1t_2 & t_1^2+2t_2
 \end{matrix} \right).$$

\vspace{0.5cm}

Let   $ \Delta_k=\det D_k$.
  \vspace{0.5cm}

  We obtain the   \textit{Infinite Different Matrix}  $D^\infty$ by continuing to operate on $d_k$ by the companion matrix,  analogous to the way we obtained $A_k^\infty$.
 \vspace{0.5cm}

 It will follow from Proposition 2 and Theorem 7 that the right hand column of  $A_k^\infty$ is the GLP sequence.

 \vspace{0.5cm}

  \begin{definition}\label{mymodp  } Let  $p$ be a rational prime,  then $p$ \textit{ramifies} if  $p$ divides $\Delta_k$.  This definition is consistent with the usual definition of ramification when speaking about number fields.

  \end{definition}

   \vspace{0.5cm}
 \begin{definition}\label{mynum.seq  } Let  $p$ be a rational prime,  let the variables $t_j$ take integer values, so that the sequences in  $\mathcal{W}$ are now numerical sequences,  and all of the matrices defined are numerical matrices.  In particular the sequences GFP and GLP become numerical sequences.  Let  $c_p$ and $c'_p$ be the periods of,  respectively, GFPmod(p) and GLPmod(p),  then we have the following theorem concerning ramification.

 \end{definition}
  \vspace{0.5cm}

   \begin{theorem}\label{myram .}  $p$ ramifies if and only if   $c_p = p \times c'_p$.

   \end{theorem}

 We shall prove this theorem in Section 6.


\section{The Convolution Product}
  \vspace{0.5cm}

   In addition to the additive structure on  $\mathcal{W}$,  there is also a multiplicative structure, namely a convolution product.  This product operates on two polynomials in $\mathcal{W}$ with the same isobaric degree to give a product of the same isobaric degree.  (In \cite{MT}  this was called the \textit{level product}.) It makes use of the fact that elements in   $\mathcal{W}$ belong to linear recursive sequences. ($\mathcal{W}$  can also be thought of as a product on sequences of weighted isobaric polynomials.)

 \vspace{0.5cm}

 \begin{definition}\label{my*prod }.Let $P_{\omega, k,n}$ and $P_{\omega', k,n}$ be two polynomials of isobaric degree n
in $\mathcal{W}$;  their convolution product is defined as
$$P_{\omega, k,n}(\mathbf{t}) \ast P_{\omega', k,n}(\mathbf{t}) = \sum_{j=0}^n  P_{\omega, k,n-j}(\mathbf{t}) P_{\omega', k,j}(\mathbf{t}).$$

\end{definition}

It is straightforward to show that this product is commutative,  associative and distributes over addition.
For the purposes of this paper this product will be used mostly to multiply a sequence by itself,  or to multiply a sequence by the sequence $(1, -t_1,\ldots,-t_j, \ldots)$. The sequence of the second option acts as an identity for this product.  We have already seen it at work in the linear recursion of Proposition 1.

If we consider the evaluations of the GFPs over the integers,  we then get a family of integer sequences (one of which is,  of course,  the Fibonacci sequence),  the product induced on these sequences turns out to be Dirichlet convolution.  Under this product the family of sequences is an abelian group, one isomorphic to the group of multiplicative arithmetic functions,  MF.   (For a proof, see \cite{MW2, LM}.)  The isomorphism was also implicitly proved in \cite{MT}  where the isobaric ring was used to construct $q$-th roots for the subgroup of rational multiplicative functions for all $q \in\mathbb{Q}$.  This turns out to be one of many ways in which the isobaric polynomials,  and in particular,  those in  $\mathcal{W}$ act as ``representing'' structures.  Of course, the Schur-hook polynomials are characters for the permutation characters of the symmetric group,  more generally,  we can represent all of the Schur polynomials in terms of the GLPs by using a version of the Jacoby-Trudi  theorem,  deriving a version of the Schur Theorem giving the character table for  $ Sym(n)$ \cite{LM}. And, number fields derived from irreducible core polynomials can also  be described completely in terms of elements in $\mathcal{W}$.  See \cite{{MW}} The first completely general analogues of the Binet formulae  were given in \cite{TM}, where again, the representation is in terms of the isobaric theory.


\section{The Isobaric Logarithm and isobaric Exponential Operators}
 \vspace{0.5cm}

 Rearick  \cite{Rear2} introduced the notions of \textit{logarithm} and \textit{exponential} with respect to the ring of arithmetic functions under the Dirichlet product.  In  \cite{LM} these notions were transported to
 the submodule of weighted isobaric polynomials in the isobaric ring.  We reproduce these definitions here:

  \vspace{0.5cm}
  $$ \mathcal{L}_n (P_{\omega,k,n}) = -t_{n-1}P_{\omega, k,1} -2t_{n-2}P_{\omega, k,2} -\cdots - t_1P_{\omega, k,n-1} +n P_{\omega, k,n} $$
 and
  $$\mathcal{E}_n(P_{\omega,k,n})= \frac{1}{n}(F_{k,n-1} P_{\omega,k,1} + F_{k,n-2} P_{\omega,k,2} + \cdots + F_{k,1} P_{\omega,k,n-1} +  P_{\omega,k,n}). $$

  It is easy to show that these operators are inverse to one another and that they have the usual properties of logarithmns and exponential functions with respect to sums and products,  in this case,  Dirichlet (convolution) products.  All of the relevant proofs appear in  \cite{LM}.  Moreover,  the following important consequence of these theorems is

    \vspace{0.5cm}

  \begin{proposition}{ \cite{LM}}\label{mypropositionlogF=G}$$\mathcal{L}_n (F_{k,n}) =  G_{k,n} $$
 and  $$\mathcal{E}_n (G_{k,n}) =  F_{k,n}.$$
  \end{proposition}

  \vspace{0.5cm}
   \begin{remark}\label{myAddF}
Proposition \ref{mypropositionlogF=G} gives the well-known result that the multiplicative group MF is isomorphic to the additive group of additive arithmetic functions.  In \cite{LM} we have called two arithmetic functions that are related by this isomorphism, \textit{Companion Sequences } (also see  \cite{DJL}).
 \end{remark}
  \vspace{0.5cm}

   Moreover, in \cite{LM} the following theorem was proved:
    \vspace{0.5cm}

    \begin{theorem}\label{mylogA=D}  $$\mathcal{L} (A_k^\infty) = (D_k^\infty).$$
    \vspace{0.5cm}
    \end{theorem}

  We can now prove the theorem in section 3:  the rational prime $p$ ramifies if and only if $c_p = p \times c'_p.$
  For, if we look at $A_k^\infty \mod(p) $ and   $D_k^\infty \mod(p)$,  then $A_k^\infty \mod(p) \cong \mathbb{Z}_{c'_p}  $ and $D_k^\infty \mod(p) \cong \mathbb{Z}_{c'_p}$. So, from the above theorem,  we have that  $\mathcal{L}$ induces the exact sequence. $$1 \rightarrow K \rightarrow \mathbb{Z}_{c_p} \rightarrow \mathbb{Z}_{c'_p} \rightarrow 1,$$
where  the kernel  $K$ is $= 0$ if and only if,  $p \in \Delta_k $.   \quad \quad   \qed

   \vspace{0.5cm}
  In \cite{LM} it was shown that the analogues to the trigonometric functions,  sine, cosine, etc. can also be defined, using the isobaric exponential function.  These isobaric trigonometric functions behave like hyperbolic trigonometric functions,  and satisfy analogues to all of the usual hyperbolic trigonometric identities.   We shall record these identities here, since we shall show in this paper that these identities can be applied to a large class of integer sequences.

   \vspace{0.5cm}
   We first define the isobaric sines and cosines (\textit{isosine, isocosine}:

     \begin{definition}\label{mylsine}
   The isosine of  $G$, $$ S(G) = \frac{1}{2}(\mathcal{E}(G) - \overline{\mathcal{E}(G)};$$
   and the isocosine of  $G$.
  $$C(G)= \frac{1}{2}(\mathcal{E}(G)  +\overline{\mathcal{E}(G)}.$$
where, $ \overline{\mathcal{E}(G)}$ means convolution inverse of $\mathcal{E}(G)$.  \end{definition}
  
    Since $\mathcal{E}(G) = F$ and $\overline{F_n }= -t_n,$  it is easy to show that

  \begin{proposition} \label{CGSG}
$$  C(G_n) =  \frac{1}{2} (F_n - {\bf t}_n);$$
$$  S(G_n) =  \frac{1}{2} (F_n +{\bf t}_n).$$
\end{proposition}\qed

Let $\delta$ be the function whose values  are $(1,0,\ldots,0,\ldots)$.

\begin{theorem}\label{mytheoremCOSSIN}
$$ C(G)^{\ast 2} - S(G)^{\ast 2} = \delta.$$
\end{theorem}
\proof \begin{eqnarray*}
&& C(G)\ast C(G)-S(G)\ast S(G) \\
& =&  \frac{1}{4}[\mathcal{E}(G)\ast \mathcal{E}(G)+\overline{\mathcal{E}(G)}\ast\overline{\mathcal{E}(G)}+2(\mathcal{E}(G)\ast\overline{\mathcal{E}(G)})]\\
&& \qquad-\frac{1}{4}[\mathcal{E}(G)\ast \mathcal{E}(G)+\overline{\mathcal{E}(G)}\ast\overline{\mathcal{E}(G)}-2(\mathcal{E}(G)\ast\overline{\mathcal{E}(G)})]\\
&=&\mathcal{E}(G)\ast\overline{\mathcal{E}(G)}\\
&=&\delta.  \end{eqnarray*}\qed
  \vspace{0.5cm}

\begin{corollary}\label{mycorollaryCOSSIN}
$$C(G_0)^{\ast 2} - S(G_0)^{\ast 2} = 1$$
$$C(G_n)^{\ast 2} = S(G_n)^{\ast 2},  n>0$$
\end{corollary}\qed

 \begin{theorem}\label{mytheoremCOS}
 Let $F$ and $G$ be induced by the core $[t_1,\ldots,t_k]$ and $F'$ and $G'$ be induced by the core $[t'_1,\ldots,t'_k]$ with $\mathcal{L}(F)=G$ and $\mathcal{L}(F')=G'$,  then
 $$C(G+G') = C(G) \ast C(G') + S(G) \ast S(G'),$$
 $$S(G+G' )= S(G) \ast C(G') + C(G) \ast S(G').$$
 \end{theorem}

\proof \begin{eqnarray*}
     & &  C(G+G')   \\
 & = &  \frac{1}{2} (\mathcal{E}(G+G') + \overline{\mathcal{E}(G+G')} ) \\
 & = &   \frac{1}{2}(\mathcal{E}(G)\ast \mathcal{E}(G') +\overline{\mathcal{E}(G)\ast \mathcal{E}(G') })
\end{eqnarray*}
while
\begin{eqnarray*}
 &  &  C(G) \ast C(G') + S(G) \ast S(G') \\
 & = &  \frac{1}{4} [(\mathcal{E}(G)+\overline{\mathcal{E}(G)})\ast(\mathcal{E}(G')+\overline{\mathcal{E}(G')}) + (\mathcal{E}(G)-\overline{\mathcal{E}(G)})\ast (\mathcal{E}(G')-\overline{\mathcal{E}(G')})] \\
& = &  \frac{1}{2}[(\mathcal{E}(G) \ast \mathcal{E}(G') + \overline{\mathcal{E}(G)} \ast \overline{\mathcal{E}(G')}]\\
& = & \frac{1}{2}[(\mathcal{E}(G) \ast \mathcal{E}(G') + \overline{\mathcal{E}(G)\ast \mathcal{E}(G')}].
\end{eqnarray*}
The proof for $S$ is analogous. \qed

\section{Multiplicative Arithmetic Functions}

We remind the reader that an arithmetic function is a function $\xi : \mathbb{N} \rightarrow \mathbb{C}$.
Such functions form a UFD under addition and Dirichlet convolution, \cite{CE}, \cite{Sha}.  An arithmetic function  $\xi$ is multiplicative if  $\xi(mn)=\xi(m)x(n) $ whenever,  $(m,n) = 1$. If the first equation holds without the second equation,  then $\xi$ is said to be \textit{completely multiplicative} CM.  The Dirichlet product on arithimetic functions is given by $$\xi_1 \ast \xi_2 (n) = \sum_{d|n} \xi_1\frac{n}{d}\xi_2 (d).$$
See, e.g., McCarthy \cite{PJM}.

 \begin{definition}\label{mydefinitionMFLR}
 An arithmetic function $f$ is \textit{locally representable}  if for each prime $p$ there is a core polynomial or power series $[t_1,t_2,\ldots,t_k,\ldots]$ such that  $f(p^n) = F_n$ for all $n\in \mathbb{Z},$  where $\{F_n\}$ is the GFP sequence induced by this core polynomial.

  \end{definition}

  \begin{definition}\label{mydefinitionMFGR}
  An arithmetic function $f$ is \textit{globally representable}  if  there is a core polynomial or power series $[t_1,t_2,\ldots,t_k,\ldots]$ such that  $f(n) = F_n$ for all $n\in \mathbb{Z},$  where $\{F_n\}$ is the GFP sequence induced by this core polynomial.

 \end{definition}

   \begin{definition}\label{mydefinitionLOCLINR}
An arithmetic function $f$ is \textit{locally linearly recursive} if for each prime  $p$
$$f(p^n) = a_1f(p^{n-1})+ a_2f(p^{n-2} )+ \cdots +  a_nf(p^0) $$
where $a_i$'s are determined by the core polynomial for  $p$.

An arithmetic function  $f$ is \textit{globally linearly recursive} if
$$f(n) = a_1f(n-1)+ a_2f(n-2) + \cdots +  a_nf(1). $$
  \end{definition}

  \begin{theorem}\label{mytheoremLOCLINR} \cite{LM}

  Multiplicative functions are locally linearly recursive.  \qed

  \end{theorem}

  \begin{remark}\label{myremarkGLR}

  If an arithmetic function is globally representable,  then it belongs to the group of units of the ring of arithmetic functions,  but lies outside of the subgroup of multiplicative arithmetic functions.  \cite{TM1,LM}.

 \end{remark}

  \begin{proposition}\label{myMFcond.}
 A necessary and sufficient condition that $f$ be representable,  either globally or locally,  is that it be, respectively, globally or locally linearly recursive.
\end{proposition}

   \proof  Since the GFP sequence is inherently a linear recursion, it is clear that the  function $f$ must also be either a locally or globally linear recursion.  On the other hand,  if it is globally or locally linearly recursive,  the parameters of the linear recursion determine a core polynomial (or power series) which in turn induces a suitable GFP sequence.  \qed
 \begin{corollary}\label{myMFLR}
  Every multiplicative function is locally linearly recursive and hence locally representable. \qed
  \end{corollary}

 \begin{proposition}\cite[Rearick, Theorem 4]{Rear1}\label{mylemmaL=0}
Let $f $ be an arithmetic function, then $f \in \mathcal{M} $ if and only if $Lf(m)  = 0$ whenever $m$ is not a power of a prime.
 \end{proposition}

 $L$ is the logarithm function for the ring of arithmetic functions defined by Rearick in the paper cited above, and $\mathcal{M}$  is the convolution group of multiplicative arithmetic functions.

\bigskip
\noindent\textbf{REMARK.}  Let  $f  \in\mathcal{A}$ and be globally and locally representable. Since $f$ is globally representable there is a global representation $f\xrightarrow{\rm gr}\,_fF$  with $\,_f\mathcal{L}(\,_fF)=\,_fG$. Since $f$ is locally representable, $f\in\mathcal{M}$ and for any prime $p$ there is a local representation $f\xrightarrow{{\rm lr}_p}\,_fF'$ with $\,_f\mathcal{L}(\,_fF')=\,_fG'$. Then
$$\,_fG_n=\left\{\begin{array}{cl}\,_fG'_r & \hbox{ if } n=p^r, \\0 & \hbox{otherwise.}\end{array}\right.$$Since $\,_fG$ and $\,_fG'$ uniquely determine $\,_fF$ and $\,_fF'$, and $\,_fF$ and $\,_fF'$ uniquely determine $f$, $f$ is well-defined and in a certain sense only trivially globally defined. So Rearick's theorem, Proposition \ref{mylemmaL=0}, tells us that an arithmetic function $f$ can not be ``non-trivially'' both locally representable and globally representable.

\section{Hessenberg Matrices}

 \begin{definition}\label{myHess  } A\textit{ lower Hessenberg matrix} is a matrix whose entries above the first super-diagonal are zero. An \textit{upper Hessenberg matrix} is defined similarly with respect to entries below the first sub-diagonal.  They are especially important for computational purposes.
\end{definition}
Our main theorem involves the following lower Hessenberg matrices:

 $$H_{+(\omega,k,n)} = \left(\begin{matrix}
t_1 & 1 & 0 & \cdots & 0 \\
t_2 & t_1 & 1 & \cdots & 0 \\
\vdots & \vdots & \vdots & \ddots & \vdots \\
t_{k-1} & t_{k-2} & t_{k-3} & \cdots & 1 \\
\omega_k t_k & \omega_{k-1} t_{k-1} & \omega_{k-2}t_{k-2}  & \cdots & \omega_1t_1
\end{matrix} \right),$$

and

$$H_{-(\omega,k,n)} = \left(\begin{matrix}
t_1 & -1 & 0 & \cdots & 0 \\
t_2 & t_1 & -1 & \cdots & 0 \\
\vdots & \vdots & \vdots & \ddots & \vdots  \\
t_{k-1} & t_{k-2} & t_{k-3} & \cdots & -1 \\
\omega_k t_k & \omega_{k-1} t_{k-1} & \omega_{k-2}t_{k-2}  & \cdots & \omega_1t_1
\end{matrix} \right),$$

Notice that these two matrices differ only in that their super-diagonals have opposite signs.

\section{Permanental and Determinantal Representations}
 \vspace{0.5cm}

 \begin{definition}\label{myperm }The \textit{permanent} of a square matrix  $M$, denoted by $perm M$ is computed by taking the determinant while ignoring the parity of the operation of the permutations on the indices.

 \end{definition}

 \vspace{0.5cm}

 \begin{theorem}\label{mypermth}
  \vspace{0.5cm}

a. $$\mbox{\rm perm }   H_{+(\omega,k,n)} = \mbox{\rm perm} \left(\begin{matrix}
t_1 & 1 & 0 & \cdots & 0 \\
t_2 & t_1 & 1 & \cdots & 0 \\
\vdots & \vdots & \vdots & \ddots & \vdots  \\
t_{k-1} & t_{k-2} & t_{k-3} & \cdots & 1 \\
\omega_k t_k & \omega_{k-1} t_{k-1} & \omega_{k-2}t_{k-2}  & \cdots & \omega_1t_1
\end{matrix} \right) = P_{\omega,k,n}(t_1,t_2, \ldots, t_k)$$
 \vspace{0.5cm}

and
 \vspace{0.5cm}

b.  $$\det  H_{-(\omega,k,n)} =\det \left(\begin{matrix}
t_1 & -1 & 0 & \cdots & 0 \\
t_2 & t_1 & -1 & \cdots & 0 \\
\vdots & \vdots & \vdots & \ddots & \vdots \\
t_{k-1} & t_{k-2} & t_{k-3} & \cdots & -1 \\
\omega_k t_k & \omega_{k-1} t_{k-1} & \omega_{k-2}t_{k-2}  & \cdots & \omega_1t_1
\end{matrix} \right) = P_{\omega,k,n}(t_1,t_2, \ldots, t_k).$$

\end{theorem}

\proof  It is convenient to follow the steps of the proof in the  $4 \times 4$- case:

$$H_{+(\omega,4,4}) = \left( \begin{matrix}
t_1 & 1 & 0 & 0 \\ t_2 & t_1 & 1 & 0 \\ t_3 & t_2 & t_1 & 1 \\ \omega_4 t_4 & \omega_3 t_3 & \omega_2 T_2 & \omega_1 t_1
\end{matrix}\right)$$
whose permanent is easily seen to be \; $\omega_1t_1^4 + (2\omega_1 +\omega_2)t_1^2t_2 + \omega_2t_2^2 + (\omega_1 + \omega_3)t_1t_3 +\omega_4t_4 = P_{\omega,4,4},$
 \vspace{0.5cm}

Moreover,  it is easy to see that there is a nesting of the Hessenberg matrices from lower right hand corner to upper left,  which enables an inductive proof of the theorem,  as the $4 \times 4$ example above suggests.
 \vspace{0.5cm}

 \vspace{0.5cm}

 Let  $M_j = \mbox{perm }H_{+(\omega,k,j)},  \, j=1,\ldots,n.$

 Clearly $M_1 = \omega_1t_1 = P_{\omega,1}$
So let  $M_j = \mbox{perm }(H_{\omega,k,j}),  j=1,\ldots ,n ,$ and assume that our result holds for $H_{+(\omega,k,j)}, j=1,\ldots, n-1,$ that is, that $M_{\omega,n-j} = P_{\omega,n-j}, j=0, \ldots,n-1.$
 \vspace{0.5cm}

\begin{lemma}\label{myM_j } $$M_n = t_1M_{n-1} + t_2M_{n-2} + \cdots + t_{n-1}M_1 + \omega_nt_n.$$
 \vspace{0.5cm}
 \end{lemma}
\proof
It is easy to see that this expression is just the computation for the permanent. But the right-hand side  $$t_1P_{\omega,n-1} + t_2P_{\omega,n-2}  + \cdots + t_{n-1}P_{\omega,1} + \omega_nt_n$$ is just   $P_{\omega,n}$,  by the linear recursion property of weighted isobaric polynomials.  So

$$P_{\omega,n} =  t_1P_{\omega,n-1} + t_2P_{\omega,n-2}  + \cdots + t_{n-1}P_{\omega,1} + \omega_nt_n = t_1M_{n-1} +  t_2M_{n-2} + \cdots +  t_{n-1}M_1 + \omega_nt_n.$$

Thus $$M_n =  t_1P_{\omega,n-1} + t_2P_{\omega,n-2}  + \cdots + t_{n-1}P_{\omega,1} + \omega_nt_n = P_{\omega,n}.$$
 
 Part b. of the theorem is proved by a similar argument. \quad \quad \qed

 The following theorem is a special case of a theorem of Kaygisiz and Sahin, \cite{KS}.

 \begin{theorem}\label{myperm=det}
 $$\mbox{perm }(H_{+(\omega,k,n)}) = \det(H_{-(\omega,k,n)}),$$
 $$\det(H_{+(\omega,k,n)}) = \mbox{perm }(H_{-(\omega,k,n)}).$$
 \end{theorem}

 Proof.  The theorem follows from parts a. and  b. of Theorem 14. \qed

 And an immediate consequence of this theorem is

  \begin{corollary}\label{myiff}
  An element of $\mathcal{W}$ has a determinantal representation if and only if it has a permanental representation.    \qed

  \end{corollary}
\section{The Representability of Integer Sequences}

The power of the isobaric theory lies in  the extent to which  other parts of mathematics are modeled0 in the isobaric ring, that is, in the ring of symmetric polynomials on the elementary symmetric polynomial basis: In \cite{MT} and \cite{MW2} It was shown that there are faithful copies of the convolution group of multiplicative arithmetic functions and of the additive group of additive arithmetic functions  in the group  $\mathcal{W}$ in the isobaric ring,  in fact, using only the GFP in the first case,  and the GLP in the second.  Moreover,  in \cite{MT}, it was shown how to embed either of these groups in their divisible closures,  using only the tools in $\mathcal{W}$ (see \cite{CG, TM}).   In \cite{MW} it was shown that all  linear recursive sequences on the integers can be represented in $\mathcal{W}$ and that the local structure of number fields  can also be described using the polynomials in $\mathcal{W}$. We have shown here that whenever an integer sequence is  linearly recursive,  or when it is the set of values of a multiplicative or an additive arithmetic function,  then it is faithfully represented by sequences in  $\mathcal{W}$.

\vspace{0.5cm}

Other examples of mathematical theories modeled in terms of $\mathcal{W}$ are given in \cite{LM}.

\vspace{0.5cm}
In this section we want to show how these results can be used to represent, and find relations among, a vast collection of integer sequences;  in fact, almost all of the well- known,  workaday sequences will appear in such a list.  This places integer sequences in a unified theory with inner structure that enables easy calculation,  and the discovery of new relations among the sequences.

\vspace{0.5cm}
It is useful at this point to give a summary of the different ways that we have of producing the elements of $\mathcal{W}$.

 \noindent1.  By using the Differential Lattice  $\mathcal{L} (t_1^{\alpha_1} t_2^{\alpha_2} \cdots t_k^{\alpha_k}).$\\
 2.  By using the permanents of $H_{+(\omega,k,n)}$.\\
 3.  By using the determanents of $H_{-(\omega,k,n)}$.\\
 4.  By  $$P_{w.k.n}(t_1,t_2,\ldots,t_k) = \omega_nt_n \ast F_{k,n} = \sum_{j=1}^n \omega_jt_jF_{k,n-j}$$\\
 5.  By $$P_{\omega,k,n}(t_1^{\alpha_1} t_2^{\alpha_2} \cdots t_k^{\alpha_k})= \sum_{\alpha \vdash n} \left( \begin{matrix}
 |\alpha|  \\  \alpha_1,\cdots,\alpha_k \end{matrix} \right)  \frac{\sum_{j=1}^k \omega_j\alpha_j}{|\alpha|}t_1^{\alpha_1}\cdots t_k^{\alpha_k}$$

 For no. 1, see \cite{MT} (remark 1 and ff.)and \cite{MT2} (p.2 ff.). It was through this lattice that the weighted isobaric polynomials were originally defined.  The differential lattice of Young diagrams is a tool frequently used by algebraic combinatorialists.Nos.  2, and 3. are described in this paper. No.  4. appeared in \cite{MT}. The fifth one is the explicit definition that results from using the methods of the differential lattice.

\section{Examples}

In this section, we first give an example of calculating the statistics of an integer sequence using the isobaric theory, followed by examples of some interesting representable integer sequences. We start with  two well-known multiplicative arithmetic functions.

1. $ \tau$ and $\sigma$ are respectively the arithmetical functions which give the number of divisors of $n$, and the sum of the divisors of $n$,  $n \in \mathbb{N}$.  They are well-known to be multiplicative functions of degree 2,  and hence are defined locally, that is, at each rational prime  $p$.

Let us see how this follows using isobaric theory:  $\tau(p^n)= n+1.$ Suppose that  $\tau$ is locally representable, that is,  suppose $\tau(p^n) = F_{k,n}(t_1,t_2,\ldots,t_k)$ for some natural number $k$ and for some integer values of $t_1,t_2, \ldots, t_k$: use the notation,  $\tau \rightarrow _{\tau}F_{k,n}$.

For simplicity in the proof we drop the initial subscript  $\tau$.  Since $p^0 = 1$,  we have $F_{k,0}=1$, and $\tau(p) = 2 = F_{k,1} = t_1$. Moreover,  $\tau(p^2) = 3 = F_{k,2} = t_1^2 + t_2$,  we infer that $t_2 = -1$.  Again,  noting that  $\tau(p^3) = 4$, we compute from the fact that $F_{k,3} = t_1^3 +2 t_1t_2 +t_3$,  that  $t_3 = 0$.  We can then compute inductively,  using that multiplicative functions are locally linearly recursive, that $t_j =0, j>2$.  Thus $k=2,$ and the core polynomial is given by $[2,-1]$.

If we apply the same line of reasoning to   $\sigma(p^n) = 1+p+\cdots+p^n$, we again deduce that $k=2$ and that the core is given by $[p+1,-p]$.  (Note that letting $p=1$ gives us back the statistics for  $\tau$ ).
The companion sequences (i.e., the isobaric logarithm), using proposition 2 and induction, are respectively,  $(2,2,\ldots,2,\ldots)$ and $(2, p+1, p^2+1, \ldots, p^n+1,\ldots)$. (Note that the value of $G_{k,0}$ is consistent with the general term, letting  $p = 1$, and, of course, with the value of the isobaric Log of $F_{k,0}$.
\vspace{0.5cm}

Using the techniques illustrated in 1. above,  we have the following interesting examples:

2.  The Euler Totient function,  $\varphi$.  $\varphi$ is MF of infinite degree; its core is a power series:  $$\varphi (p^n) = p^n - p^{n-1}$$ $$ [t_1,\ldots, t_k,\ldots] = [p-1, p-l,\ldots, p-1,\ldots].$$

It is locally $F$-representable.  Its companion sequence,  which is locally $G$-representable, is given by
$$G_0 = 2, \quad  G_n = p^n-1, \quad  n>0. $$ The trigonometric identities then apply.  We leave these computations to the interested reader.

3.  Jordan's function,  $J_k(n)$ is a multiplicative function and a generalization of the Euler function:  $$J_k(n)= p^{kn}-p^{(k-1)n}.$$   It's degree is infinite, that is,  $[p^k-1,\ldots, p^k-1,\ldots].$, and it is locally $F$-representable.  Its companion sequence is given by $$G_0 = 2,\quad G_n = p^{kn}-1,\quad  n>0.$$

4.  The Catalan sequence  $C_n = \frac{1}{n+1}\left( \begin{matrix}
2n  \\ n+1 \end{matrix} \right)$ is  $F$-globally representable with core $[t_j = C_{j-1}]$. We have named such a sequence, \textit{incestuous}.  It's companion sequence is number A001700in the Encyclopedia  of Integer Sequences Sequences and is given by Proposition 2 as
$$_CG_n =\left(\begin{array}{c }
     2n-1    \\
      n
\end{array}\right).$$
\vspace{0.5cm}

It also has interpretations in terms of Dyck paths.
\vspace{0.5cm}

5.  The Chebychev polynomials  of the 2nd kind, $U_n(x)$,  are globally $F$-representable with core $[2x,-1]$. If, however, we let $x=1$,  then we get the local $F$-representation of $\tau$..  Multiplicative functions belong to the group of units of the ring of arithmetic functions.  However,  the group of units is much larger than the group of multiplicative functions  \cite{TM1} .  In particular, any  $F$-globally representable function belongs to this group.  Hence, for any value of $x$,  $\{U_n(x)\}$ we obtain an element of this group.
\vspace{0.5cm}

6.   The Stirling numbers of the 2nd kind, usually denoted as
$$\left\{\begin{array}{c }
     n    \\
      k
\end{array}\right\},$$
for which a kind of ``Pascal'' recursion relation holds; namely,
$$\left\{\begin{array}{c }
     n+1    \\
      k
\end{array}\right\}=k\left\{\begin{array}{c }
     n    \\
      k
\end{array}\right\}+\left\{\begin{array}{c }
     n    \\
      k-1
\end{array}\right\}.$$

$\begin{array}{ cccccccccccc}
{\color{red}     n/k }  &   0  & 1 & 2 & 3 & 4 & 5 & 6 & 7 & 8 & 9 & 10\\
  {\color{red}    0}    &  1   \\
 {\color{red}     1}    &  0  &  1 \\
 {\color{red}     2}    &  0  &  1  &  1 \\
  {\color{red}    3}    &  0  &  1  &  3  &  1 \\
 {\color{red}     4}    &  0  &  1  &  7  &  6  &  1 \\
 {\color{red}     5}    &  0  &  1  &  15 & 25  &  10  &  1  \\
    {\color{red}  6}    &  0  &  1  &  31 &  90  &  65  &  15  &  1  \\
    {\color{red}  7}    &  0  &  1  &  63  &  301 &  350 & 140 & 21 & 1 \\
    {\color{red}  8}    &  0  &  1  & 127 & 966  & 1701 & 1050 & 266 & 28 & 1 \\
    {\color{red}  9}    &  0  &  1  &  255 & 3025 &7770 & 6951 & 2646 & 462 & 36 & 1 \\
     {\color{red} 10}  &  0  &  1  &  511 & 9330 & 34105 & 42525 & 22827 & 5880 & 750 & 45 & 1 \\
\end{array} $

\vspace{0.5cm}

 (Table of Stirling Numbers of the second kind, $n=0,\ldots,10$)

 \vspace{0.5cm}
Note that the first non-zero element in each column is a 1,  and that the set of next-to-last numbers in each row gives the triangular numbers in the correct ascending sequence.

   \begin{theorem}\label{myStirling2}
 \vspace{0.5cm}

   Each column is an an arithmetic function in the group of units of the ring of arithmetic functions; that is each column is a globally $F$-representable, degree $k$,  arithmetic function.  For example, the first few core polynomials are given by
   $$[1], [3,-2], [6,-11,6], [10,-35, 50, -24], [15,-85, 225,-274,120] $$
One will recognize these numbers as the coefficients of the polynomials  otherwise given as $$(x-1);  (x-1)(x-2):  (x-1)(x-2)(x-3); (x-1)(x-2)(x-3)(x-4);  (x-1)(x-2)(x-3)(x-4)(x-5), $$ This is not a surprise:  the Stirling numbers of the second kind are given as the coefficients of the polynomials defined by so-called ``falling'' factorial.
   That is  $$\mathcal{C}(X) = \prod_{j=1}^k (X-j).$$
 $t_{k,1} = T_k$ where $t_{k,1}$ is just the recursion parameter  $t_1$ for the $k$-th column,  and  $T_k $ is the $k$-th element in the sequence of triangular numbers.  The last coefficient in the core representation  is $t_k = k!$.
\end{theorem}

As for Stirling numbers of the first kind, they are not representable in our sense,  however,  there is an unexpected and interesting relation between Stirling numbers of the 1st kind and Stirling numbers of the 2nd kind that we do not understand.  Namely,  the rows,  in what might be called the Stirling triangle of the first kind,  are just the absolute values of the core numbers,  $t_j$  in reverse order. For example,  the 5th row of the table for Stirling  numbers of the 1st kind is just $$(0, 24,50, 35, 10, 1)$$ (The apparently extra one is just $t_0$,  the leading coefficient of the core polynomial,  which we have suppressed in our notation).

As a result of these observations,  we have the following relation between Stirling numbers of the first kind and Stirling numbers of the second kind.

   \begin{corollary}\label{myStirling1&2}

   $$\sum_{j=0}^{k-1} (-1)^{j-1}  \left\{\begin{array}{c }
     n+j    \\
      k
\end{array}\right\}      \left[\begin{array}{c }
     k+1    \\
      k-j
\end{array}\right]  = \left\{\begin{array}{c }
     n+k    \\
      k
\end{array}\right\},$$
where we use   $   \left[\begin{array}{c }
     k+1    \\
      k-j
\end{array}\right] $ to denote Stirling numbers of the first kind.

   \end{corollary}

   That is,  the recursion parameters of the $k$-th column of the table of Stirling numbers of the 2nd kind is  $$t_{k-1} = (-1)^{j-1}   \left[\begin{array}{c }
     k+1    \\
      k-j
\end{array}\right]  . $$

\proof  The proof follows from the definitions of the two types of sequences in terms of rising and falling factorials. \qed
\vspace{0.5cm}

The companion sequence for the $k$-th column is just the sequence $1^n+2^n+\cdots + k^n$.  When  $k=3$, for example, this is the sequence A001550 in \cite{EIS} .
   \vspace{0.5cm}

7. The Trigonometric series \textsl{$T_k$}  is also $F$-globally  representable:    $\mathcal{C}(X) = [3,-2,1]$, so that  $k=3$.  The companion sequence is just  $\{G_n = 3\}$ for all $n$.

\vspace{0.5cm}

In \cite{KS},  it was pointed out that Pell numbers, Pell-Lucas numbers, bivarate Fibonacci numbers, Perrin sequences, and Exponential Perrin sequences have permanental and determinantal representations.  We give their statistics here.
\vspace{0.5cm}

8.   Pell numbers are $F$-representable with a core $[2,1]$. \\

9.   Pell-Lucas numbers are their companions,  hence, $G$-representable. \\

10. Bivarate Fibonacci numbers, are $F$-representable with  $k=2$. \\

11. Perrin sequences are $G$-representable with a core of  $[0,1,1]$. \\

12. Exponential Perrin sequences are $F$-representable on the same core,  thus 11. is the companion sequence for 12.


\section{Apendix: tables of the first few GFP and GLP}

In these tables,  we take the point-of-view that $k$ goes to $\infty$,  and we omit it from the indices of the polynomials. When the companion matrix is non-singular,  all sequence can be extended "nort5hward",  this extension will also be omit from the tables below.

\vspace{0.5cm}

 \flushleft
 $F_0 = 1 $ \\
$F_1 = t_1 $\\
$F_2 = t_1^2 + t_2$ \\
$F_3 = t_1^3 + 3t_1t_2 + t_3 $\\
$F_4 = t_1^4 + 3t_1^2t_2 + t_2^2 + 2t_1t_3 + t_4$ \\
$F_5 = t_1^5 + 4t_1^3t_2 + 2t_1t_2^2 + 2t_2t_3 + 2t_1t_4 + t_5 $\\
$F_6 = t_1^6 + 5t_1^4t_2 + 6t_1^2t_2^2 + t_2^3 + 4t_1^3t_3 + t_3^2 + 6t_1t_2t_3 + 3t_1^2t_4 + 2t_2t_4 + 2t_1^2t_5 + t_6 $ \\
\qquad$\vdots$
\vspace{0.5cm}

$G_0 = k$ \\
$G_1 = t_1$ \\
$G_2 = t_1^2 + 2t_2$ \\
$G_3 = t_1^3 + 3t_1t_2 + 3t_3$ \\
$G_4 = t_1^4 + 4t_1^2t_2 +2 t_2^2 + 4t_1t_3 + 4t_4$ \\
$G_5 = t_1*5 + 5t_1^5t_2 + 5t_1t_2^2 + 5t_2t_3 + 5t_1t_4 + 5 t_5 $\\
$G_6 =  t_1^6 + 6t_1^4t_2 + 6t_1^2t_2^2 + 2t_2^3 + 6t_1^3t_3 +3 t_3^2 + 12t_1t_2t_3 + 6t_1^2t_4 + 6t_2t_4 + 6t_1^2t_5 + 6t_6 $ \\
\qquad$\vdots$

Huilan Li
Department of Mathematics, Drexel University, Philadelphia, PA 19104 U.S.A.
Email: \textit{huilan.li@gmail.com}

Trueman MacHenry
Department of Mathematics and Statistics, York University, Toronto, Ontario M3J 1P3 CANADA
machenry@mathstat.yorku.ca	
 \end{document}